\documentclass[12pt,leqno]{article}
 \usepackage{amsmath,amssymb,amsthm}
\newcommand{\red}{\protect\large\bf}

\newtheorem{thm}{Theorem}

\newtheorem{lem}[thm]{Lemma}

\begin{document}
 \font\klsl=cmsl10 at 9pt
\font\fraktur=eufm10

\title{\textrm Contact 3-manifolds and Ricci solitons}
\author{Jong Taek Cho}
\date{}
\maketitle
 \footnote{\textit{Date}: May 7th, 2010}
%\footnote{}
%\keywords{a} \subjclass{30}

\begin{abstract}
\noindent
 A contact 3-manifold $M$ admitting a transversal Ricci soliton $(g,v,\lambda)$ is either Sasakian  or  locally isometric to one of the Lie groups $SU(2)$, $SL(2,R)$, $E(2)$, $E(1,1)$ with a left invariant metric.
\end{abstract}
\section{\red Introduction}
Contact geometry is motivated by classical mechanics. Where a symplectic space is considered  the even-dimensional phase space of a mechanical system, a contact space corresponds to the odd-dimensional extended phase space that includes the time variable. A {\it contact manifold} $(M,\eta)$ is a smooth manifold $M^{2n+1}$ together with a global one-form $\eta$ such that $d\eta$ has maximal rank $2n$ on the contact distribution $D=\hbox{ker }\eta$. The duality of $\eta$ defines a unique vector field $\xi$, the {\it Reeb vector field}. The {\it Reeb flow} is a one parameter group of diffeomorphisms $\{\phi_t\}$ generated by the Reeb vector field $\xi$.
 Martinet \cite{Martinet} proved that every closed, orientable 3-manifold admits a contact structure.
\vskip.15cm
A \textit{Ricci soliton} is defined on a Riemannian manifold ($M,g$) by
\begin{equation}\label{Ricci}
\frac{1}{2}\pounds_W g +  \operatorname{Ric} - \lambda g = 0
\end{equation}
where $W$ is a vector field (the potential  vector field), $\lambda$ a constant on $M$. A Ricci soliton with $W$ zero is reduced to Einstein equation. Not only for studying the topology of manifolds, but in the study of string
theory, theoretical physicists also have been looking into the equation of Ricci solitons
\cite{Friedan}.
Compact Ricci solitons are the fixed points of the \textit{Ricci flow}:
\begin{equation}\label{Ricciflow}
\displaystyle\frac{\partial}{\partial t}g = -2\operatorname{Ric}
\end{equation}
 projected from the space of metrics onto its quotient modulo diffeomorphisms and scalings, and often arise as blow-up limits for the Ricci flow on compact manifolds. Indeed, the Ricci flow $\{g_t\}$ is equivalent to the initial metric $g = g_0$ satisfying the Ricci soliton equation (1) for some vector field $W$ and some constant $\lambda$. The Ricci soliton is said to be shrinking, steady, and expanding according as $\lambda>0$, $\lambda=0$, and $\lambda<0$ respectively.  For details we refer to \cite{Cao} or \cite{ChowKnopf} about the Ricci flows and their solitons. Hamilton \cite{Hamilton} initiated the Ricci flow theory, namely, he proved that on a compact manifold $M$, the Ricci flow equation with
any prescribed initial metric $g_0$ has a unique solution on some maximal time
interval $[0, T)$, where $0 < T\leq \infty$. Developing it, Hamilton proved that
 a compact 3-manifold admitting a metric of strictly positive Ricci
curvature admits in fact a metric of constant positive sectional curvature. In particular, if such a manifold is simply connected, the manifold is isometric to the sphere, and then is diffeomorphic to the $S^3$. Then the Ricci flow has been a very attractive approach to a possible
positive answer to the Poincar\'e Conjecture \cite{Poincare}, more generally to the Thurston's Geometrization Conjecture \cite{Thurston}, which describes all 3-dimensional manifolds in terms of the Eight Geometries (cf. Remark 2).

% After all, implementing the program Perelman \cite{Perel1}-\cite{Perel3} proves the Poincar\'e conjecture.

 %Hamilton \cite{Hamilton} and Ivey \cite{Ivey} proved that a Ricci soliton on a compact manifold has constant curvature in dimension 2 and 3, respectively. If the vector field $V$ is the gradient of a potential function, then $g$ is called a gradient Ricci soliton. We refer to \cite{Cao} and \cite{ChowKnopf} for details about Ricci solitons or gradient Ricci solitons.
 % and equation (\ref{Ricci}) assumes the form
%\begin{equation}\label{GradRicci}
%\nabla\nabla f + \operatorname{Ric} - \lambda g =0
%\end{equation} We refer to \cite{ChowKnopf} for details about Ricci solitons or gradient Ricci solitons.
%In \cite{Perel}(Remark 3.2) Perelman showed that a Ricci soliton on a compact manifold is always %\textit{replaced by} a gradient Ricci soliton.
%Indeed, on a compact Riemannian manifold  $(M,g)$ a potential vector field $V$  is the sum of
%a Killing field and a gradient %Hence, (\ref{Ricci})
%becomes  $g(\nabla_X Df,Y) +  \operatorname{Ric(X,Y)} =\lambda g(X,Y)$  for  a function $f$, where   $Df$ %denotes the gradient of $f$.
%(see also, \cite{Derdzinski}, \cite{LopezRio}).
% and equation (1) assumes the form
%\begin{equation*}
%\nabla\nabla f = \operatorname{Ric}- \lambda g
%\end{equation*}
\vskip.15cm
Now, we consider the Ricci flow in contact geometry \cite{Cho}. Then, we keeping Martinet's result in mind,
it is very natural and interesting to study the Ricci flow in contact 3-manifolds. In this context, we first establish a very special Ricci flow which evolves by the Reeb flow and a (time dependent) evolving factor at the same time. Then we have the corresponding Ricci soliton equation, which we call a {\em contact Ricci soliton}:
\begin{equation}\label{contRicci}
\frac{1}{2}\pounds_\xi g +  \operatorname{Ric} - \lambda g = 0.
\end{equation}
In \cite{Cho}, we proved
\noindent
\begin{thm} A 3-dimensional contact Ricci soliton $(g,\xi,\lambda)$ is of constant curvature $+1$.
\end{thm}
Next, as a complementary partner of the contact Ricci soliton,  we consider the so-called {\it transversal Ricci soliton}:
\begin{equation}\label{contRicci1}
\frac{1}{2}\pounds_v g +  \operatorname{Ric} - \lambda g = 0,
\end{equation}
where $v$ is a complete vector field orthogonal to $\xi$.
Then, we prove in Section 3 the following theorem.

\noindent
\begin{thm} A contact 3-manifold $M$ admitting a transversal Ricci soliton $(g,v,\lambda)$ is either Sasakian or locally isometric to one of the following Lie groups with a left invariant metric: $SU(2)$, $SL(2,R)$, $E(2)$$($the group of rigid motions of the Euclidean 2-space$)$, $E(1,1)$$($the group of rigid motions of the Minkowski 2-space$)$.
\end{thm}
\vskip.15cm
%For a given contact structure $\eta$ in a closed 3-manifold $M$, it is known that there is a compatible Riemannian metric $g$.
%Applying the Hamilton's Ricci flow, then together with Theorem 1 we can reach at the Poincar\'e conjecture:
%\vskip.2cm
%\textit{Every closed, smooth, simply connected 3-manifold is diffeomorphic to a sphere.}
% Hence, we have
%\vskip.2cm
%\noindent
%{\bf Corollary B}. {\it A contact Ricci soliton is a shrinking soliton.}
%\vskip.15cm

 \section{\red Contact geometry and the known results}

We start by reviewing briefly the fundamental materials about contact
geometry. All manifolds in the present paper are assumed to
be connected, oriented and smooth.
 \vskip.2cm
 A $(2n+1)$-dimensional manifold $M$ is a {\it contact manifold} if it is equipped with a global one-form $\eta$ such that
 $\eta\wedge (d\eta)^n \neq 0$ everywhere.
 Given a contact form $\eta$, there
exists a unique vector field $\xi$, called the {\it Reeb
vector field}, satisfying $\eta(\xi)=1$ and $d\eta(\xi,X)=0$ for
any vector field $X$. It is well-known that there also exists a
Riemannian metric $g$ and a $(1,1)$-tensor field $\varphi$ such
that
\begin{equation}\label{2.1}
\eta(X)=g(X,\xi),\ d\eta(X,Y)=g(X,\varphi Y),
        \ \varphi^2 X=-X+\eta(X)\xi,
\end{equation}
where $X$ and $Y$ are vector fields on $M$. From (\ref{2.1}), it
follows that $\varphi\xi=0,\ \eta\circ\varphi=0,
        \ g(\varphi X,\varphi Y)=g(X,Y)-\eta(X)\eta(Y).$
A Riemannian manifold $M$ equipped with structure tensors $(\eta,g)$
satisfying (\ref{2.1}) is said to be a {\em contact Riemannian
manifold} or a {\em contact metric manifold} and it is denoted by
$M=(M,\eta,g)$. Given a contact Riemannian manifold $M$, we define a
$(1,1)$-tensor field $h$ by $h=\frac12 \pounds_{\xi}\varphi$. Then $h$ is self-adjoint and
satisfies
\begin{equation}\label{2.3}
                 h\xi=0,\quad  h\varphi=-\varphi h,
\end{equation}
\begin{equation}\label{2.4}
                 \nabla_X\xi=-\varphi X-AX,
\end{equation}
 where $\nabla$ is Levi-Civita connection and $A=\varphi h$.
From (\ref{2.3}) and (\ref{2.4}), we see that each trajectory of $\xi$ is a
geodesic flow.
%Furthermore, we know that $\nabla_\xi \varphi=0$ in
%general (cf. p. 67 in \cite{Blair}). From the second equation of
%(\ref{2.3}) it follows also that
%\begin{equation}
%\label{2.4-1}(\nabla_\xi h)\varphi=-\varphi (\nabla_\xi
%h).
%\end{equation}
% We denote by $R$ the Riemannian curvature tensor
%defined by
%$$
%   R(X,Y)Z=\nabla_X(\nabla_Y Z)-\nabla_Y(\nabla_X Z)-\nabla_{[X,Y]}Z
%$$
%for all vector fields $X,Y,Z$. Along a characteristic flow $\xi$,
%the Jacobi operator $\ell=R(\cdot,\xi)\xi$ is a symmetric
%$(1,1)$-tensor field. We call it the {\em characteristic Jacobi
%operator}. From the definition of $R$ by using (\ref{2.4}) we have
%\begin{equation}
%  \label{2.5} \ell=-\varphi^2+\varphi\nabla_\xi h-h^2.
%\end{equation}
%From (\ref{2.5}) using the 2nd equation of (\ref{2.3}) and
%(\ref{2.4-1}) we have
%\begin{equation}
%  \label{2.6} \nabla_\xi h=1/2(\ell\varphi-\varphi \ell).
%\end{equation}
Moreover, we also have
\begin{equation}\label{2.5}
\operatorname{Ric}(\xi,\xi)=2n- \hbox{trace } h^2
\end{equation}
(cf. Corollary 7.1 in \cite{Blair}).
A contact Riemannian manifold for which $\xi$ is Killing is called a
{\em $K$-contact manifold.} It is easy to see that a contact
Riemannian manifold is $K$-contact if and only if $h=0$.
For a contact manifold $M$, the tangent space $T_pM$ of $M$ at each
point $p\in M$ is decomposed as $T_pM=D_p\oplus\{\xi\}_p$(direct
sum), where we denote $D_p=\{v\in T_pM|\eta(v)=0\}$. Then the
$2n$-dimensional distribution (or subbundle) $D:p\rightarrow D_p$ is
called the {\em contact distribution $($or contact subbundle$)$}.
For a contact manifold $M$, the associated almost CR structure is given by the holomorphic
subbundle ${\cal H}=\{ X-i JX: X\in D \}$ of the
complexification ${TM}^C$ of the tangent bundle $TM$, where
$J=\varphi|D$, the restriction of $\varphi$ to $D$. We say that {\em the almost CR structure is integrable} if
$[{\cal H},{\cal H}]\subset {\cal H}$. A {\em Sasakian manifold is a $K$-contact manifold whose associated almost CR structure is integrable}. Then we observe that a 3-dimensional $K$-contact
 manifold is already Sasakian.
\vskip.15cm

 Recall another class of contact metric manifolds, the so-called {\em contact $(\alpha,\beta)$-manifolds} (introduced by Blair, Koufogiorgos and Papantoniou \cite{BKP}) are defined by the curvature condition
\begin{equation*}
R(X,Y)\xi=\alpha\big(\eta(Y)X-\eta(X)Y\big)+\beta\big(\eta(Y)hX-\eta(X)hY\big)
\end{equation*}
for arbitrary vector fields $X,Y$ and for some real numbers $\alpha$ and $\beta$.
The class of contact ($\alpha,\beta$)-manifolds is developed from a contact metric manifold with $R(X,Y)\xi=0$.
Indeed, by a $D$-homothetic deformation of a contact metric manifold with $R(X,Y)\xi=0$, we obtain a contact ($\alpha,\beta$)-manifold. Here, for a positive constant $\epsilon$, $D_{\epsilon}$-homothetic deformation means a change of structure tensors by
$$
\bar{\eta}=\epsilon\eta,\quad  \bar{\xi}=\frac{1}{\epsilon}\xi,\quad  \bar\varphi=\varphi,\quad  \bar{g}=\epsilon g+\epsilon(\epsilon-1)\eta\otimes\eta.
$$
This class includes Sasakian manifolds (for $\alpha = 1$ and $h=0$) and the trivial sphere bundle $E^{n+1}\times S^{n}$ (for $\alpha=\beta=0$). Characteristic examples of non-Sasakian ($\alpha,\beta$)-contact spaces are the tangent sphere bundles of Riemannian manifolds of constant curvature $\ne 1$.
For the three dimensions, such spaces are classified in \cite{BKP}.
\begin{thm}\label{BKP}
A 3-dimensional contact $(\alpha,\beta)$-space is either Sasakian or  locally isometric to one of the Lie groups  $SU(2)$, $SL(2,R)$, $E(2)$, $E(1,1)$ with a left invariant metric.
\end{thm}

Boeckx and the present author \cite{BC} proved the following theorem, which has a crucial role in proving our main Theorem 2.

\begin{thm}\label{BC}
A 3-dimensional contact manifold with $\eta$-parallel $h$, which means $g((\nabla_x h)y,z)=0$ for any vector fields $x,y,z$ orthogonal to $\xi$, is a contact $(\alpha,\beta)$-space.
\end{thm}
We refer to \cite{Blair} for the above formulas, results and the further details on contact Riemannian geometry.
\section{\red Proof of Theorem 2}

%First of all, we prepare a general formula for Ricci solitons.
%\begin{lem}$(\cite{Cho})$
%If $(g,V)$ is a Ricci soliton of a Riemannian manifold, then we have
%\begin{equation}\label{norm}
%\frac{1}{2}\|\pounds_V g\|^2=dr(V)+2\operatorname{div}(\lambda V-SV),
%\end{equation}
%where $r$ denotes the scalar curvature of $g$ and $S$ the Ricci operator defined by $\operatorname{Ric}(X,Y)=g(SX,Y)$.
%\end{lem}
%\begin{proof}
%We adapt a local coordinate system $(x^i)$. Then the equation (\ref{Ricci}) entails
%\begin{equation}\label{tail}
%\frac{1}{2}\pounds_V g^{ij} +  R^{ij} - \lambda g^{ij} = 0.
%\end{equation}
%From the above equation (\ref{tail}) we compute
%\begin{equation}\label{1}
%\begin{split}
%\frac{1}{2}\|\pounds_V g\|^2&=- R^{ij}\pounds_V g_{ij} + \lambda g^{ij}\pounds_V g_{ij}\\
%&=-\pounds_V r+g_{ij}\pounds_V R^{ij}+\lambda g^{ij}\pounds_V g_{ij}.\\
%&=-\pounds_V r+g^{ij}(\nabla_kV_i{R^k}_{j}+\nabla_kV_j{R^i}_{k})-\lambda g^{ij}(\nabla_kV_i{g^k}_{j}+\nabla_kV_j{g^i}_{k}).
%\end{split}
%\end{equation}
%We compute the second term of the last equation:
%\begin{equation}\label{21}
%\begin{split}
%g_{ij}\pounds_V R^{ij}
%&=g_{ij}\nabla_V R^{ij}-g_{ij}\nabla_{\alpha}V^i R^{\alpha j}-g_{ij}\nabla_{\alpha}V^j R^{i\alpha}\\
%&=g_{ij}\nabla_V R^{ij}-2g_{ij}\nabla_{\alpha}V^i R^{\alpha j}\\
%&=2dr(V)-2\operatorname{div}{SV},
%\end{split}
%\end{equation}
%where we have used $dr(V)=2V^\beta\nabla_\alpha R^{\alpha}_{\ \beta}$.
%Since $g^{ij}\pounds_V g_{ij}=2\operatorname{div}{V}$, using (\ref{1}) and (\ref{21}) we obtain (\ref{norm}).\end{proof}

 We consider on $M$ the maximal open subset $\mathcal{U}_1$ on which $h\neq0$ and the maximal open subset $\mathcal{U}_2$ on which $h$ is identically zero. Then $\mathcal{U}_1\bigcup\mathcal{U}_2$ is open dense in $M$. Suppose that $M$ is non-Sasakian. Then $\mathcal{U}_1$ is non-empty and there exists a local orthonormal frame field
$\mathcal{E}=\{e_1,e_2=\varphi e_1,e_3=\xi\}$ such that
$
he_1=\mu e_1,\ he_2=-\mu e_2,
$
where $\mu$ is a smooth function.
First of all, we prepare the following lemma.
\begin{lem}\label{LemmaGX}
{\rm(\textit{cf.} \cite{CPV})}
Let $M$ be a $3$-dimensional contact Riemannian
manifold.
Then with respect to $\mathcal{E}$, the Levi-Civita
connection $\nabla$ is given by
$$
\begin{array}{ccc}
\nabla_{e_1}{e_1}= be_{2},
&
\nabla_{e_1}{e_2}=-be_{1}+(1+\mu)\xi,
&
\nabla_{e_1}\xi=-(1+\mu)e_2,\\
\nabla_{e_2}{e_1}= -ce_{2}+(\mu-1)e_{3},
&
\nabla_{e_2}{e_2}=ce_{1},
&
\nabla_{e_2}\xi=(1-\mu)e_1,\\
\nabla_{\xi}{e_1}=a e_2,
&
\nabla_{\xi}{e_2}=-a e_1,
&
\nabla_{\xi}\xi=0,
\end{array}
$$
where $a,b,c$ are smooth functions.
The Ricci operator $S$ is given by
\begin{eqnarray*}
Se_1&=& \operatorname{Ric}(e_1,e_1)e_1+\xi(\mu)e_{2}+
(2b\mu-e_{2}(\mu))\xi,
\\
Se_{2}&=&
\xi(\mu)e_{1}+
\operatorname{Ric}(e_2,e_2)e_{2}+
(2c\mu-e_{1}(\mu))\xi,\\
S\xi&=&
(2b\mu-e_{2}(\mu))e_{1}+
(2c\mu-e_{1}(\mu))e_{2}+
2(1-\mu^{2})\xi.
\end{eqnarray*}
%and
%$$\operatorname{Ric}(e_1,e_1)=\frac{r}{2}+\mu^{2}
%-2a\mu-1,\ \
%\operatorname{Ric}(e_2,e_2)=\frac{r}{2}+\mu^{2}
%+2a\mu-1,
%$$
%where $r$ denotes the scalar curvature of $M$.
\end{lem}
%\vskip.2cm
%The following lemma has a crucial role in proving Theorem 2.
%\begin{lem}$(\cite{Cho})$
%If $(g,W)$ is a Ricci soliton of a Riemannian manifold, then we have
%\begin{equation}\label{norm}
%\frac{1}{2}\|\pounds_W g\|^2=dr(X)+2\operatorname{div}(\lambda X-SX),
%\end{equation}
%where $r$ denotes the scalar curvature of $g$ and $S$ the Ricci operator defined by $\operatorname{Ric}(X,Y)=g(SX,Y)$.
%\end{lem}
%\begin{proof}
%We adapt a local coordinate system $(x^i)$. Then the equation (\ref{Ricci}) entails
%\begin{equation}\label{tail}
%\frac{1}{2}\pounds_W g^{ij} +  R^{ij} - \lambda g^{ij} = 0.
%\end{equation}
%From the above equation (\ref{tail}) we compute
%\begin{equation}\label{1}
%\begin{split}
%\frac{1}{2}\|\pounds_W g\|^2&=- R^{ij}\pounds_W g_{ij} + \lambda g^{ij}\pounds_W g_{ij}\\
%&=-\pounds_W r+g_{ij}\pounds_W R^{ij}+\lambda g_{ij}\pounds_W g^{ij}.\\
%&=-\pounds_V r+g^{ij}(\nabla_kV_i{R^k}_{j}+\nabla_kV_j{R^i}_{k})-\lambda g^{ij}(\nabla_kV_i{g^k}_{j}+\nabla_kV_j{g^i}_{k}).
%\end{split}
%\end{equation}
%We compute the second term of the last equation:
%\begin{equation}\label{2}
%\begin{split}
%g_{ij}\pounds_W R^{ij}
%&=g_{ij}\nabla_W R^{ij}-g_{ij}\nabla_{\alpha}W^i R^{\alpha j}-g_{ij}\nabla_{\alpha}W^j R^{i\alpha}\\
%&=g_{ij}\nabla_W R^{ij}-2g_{ij}\nabla_{\alpha}W^i R^{\alpha j}\\
%&=2dr(W)-2\operatorname{div}{SW},
%\end{split}
%\end{equation}
%where we have used $dr(W)=2W^\beta\nabla_\alpha R^{\alpha}_{\ \beta}$.
%Since $g_{ij}\pounds_W g^{ij}=2\operatorname{div}{W}$, using (\ref{1}) and (\ref{2}) we obtain (\ref{norm}).
%\end{proof}

\vskip.15cm
%A contact Riemannian manifold $(M,\eta,g)$ is said to be $\eta$-Einstein if its Ricci tensor is of the form
%$$\operatorname{Ric}=\alpha g+\beta \eta\otimes\eta$$ for functions $\alpha,\beta$ on $M$. For an $\eta$-Einstein $K$-contact manifold in dimension $>3$,
%it is known that $\alpha,\beta$ are constants, and hence the scalar curvature is constant

For $\mathcal{E}$, it is known that $e_3=\xi$ is defined globally on $M$. Lifting to the universal covering space $\widetilde{M}^3$ if necessary  we have global orthonormal frame field, which also denoted by $e_1$, $e_2$ and $e_3$.
Now, we assume that ${M}^3$ admits a transversal Ricci soliton $(g,v,\lambda)$ with $v=f_1 e_1+f_2 e_2$, where $f_1$ and $f_2$ are smooth functions.
Then the Ricci soliton equation (\ref{contRicci1}) is written by
\begin{equation}\label{Ricci1}
\frac{1}{2}\Big(g(\nabla_X v,Y)+g(\nabla_Y v,X)\Big) +  \operatorname{Ric}(X,Y) - \lambda g(X,Y) = 0.
\end{equation}
We compute
$$\nabla_X v=(Xf_1)e_1+f_1\nabla_X e_1+(Xf_2)e_2+f_2 \nabla_X e_2.$$
Putting $X=Y=\xi$ in (\ref{Ricci1}), then by using (\ref{2.5}) we get $\lambda=2-2\mu^2,$ from which we find that $\mu$ is constant.
Put $X=Y=e_1$ to get
\begin{equation}\label{2.7}
e_1(f_1)-bf_2+\operatorname{Ric}(e_1,e_1)=\lambda.
\end{equation}
Putting $X=Y=e_2$, then we obtain
\begin{equation}\label{2.8}
e_2(f_2)-cf_1+\operatorname{Ric}(e_2,e_2)=\lambda.
\end{equation}
If we put $X=\xi$ and $Y=e_1$ in (\ref{Ricci1}), then using the formula in Lemma \ref{LemmaGX} we get
\begin{equation}\label{2.6}
\frac{1}{2}\Big(\xi(f_1)-a f_2+(1+\mu)f_2\Big)+\operatorname{Ric}(\xi,e_1)=0.
\end{equation}
Put $X=\xi$ and $Y=e_2$ in (\ref{Ricci1}) to get
\begin{equation}\label{2}
\frac{1}{2}\Big(a f_1+\xi(f_2)+(\mu-1)f_1\Big)+\operatorname{Ric}(\xi,e_2)=0.
\end{equation}
Since $\operatorname{Ric}(e_1,e_2)=0$, if we put $X=e_1$ and $Y=e_2$ in (\ref{Ricci1}), then we  obtain
\begin{equation}\label{2.10}
bf_1+e_1(f_2)+e_2(f_1)+cf_2=0.
\end{equation}
\vskip.15cm
 We suppose that $f_1$ and $f_2$ are pointwise linearly independent functions. Give an additional assumption:
\begin{equation}\label{2.11}
e_1(f_2)+e_2(f_1)=0.
\end{equation}
Then, from (\ref{2.10}) we have $b=c=0$. From the computation $(\nabla_{e_i}h) e_j=\nabla_{e_i} (h e_j)-h\nabla_{e_i}e_j$ for $i,j=1,2$, we can easily
obtain $g((\nabla_{e_i}h) e_j,e_j)=0$, $g((\nabla_{e_1}h) e_1,e_2)=2\mu g(\nabla_{e_1} e_1, e_2)=0$ and $g((\nabla_{e_2}h) e_1,e_2)=2\mu g(\nabla_{e_2} e_1, e_2)=0$,
where we have used $\mu$ is constant and $b=c=0$. Namely, we have
\begin{lem}\label{Lemma6}
 $h$ is $\eta$-parallel, i.e. $g((\nabla_x h)y,z)=0$ for any vector fields $x,y,z$ orthogonal to $\xi$.
\end{lem}
Hence, from Theorem \ref{BC} we have $M$ is a contact $(\alpha,\beta)$-space. Then, due to Theorem \ref{BKP} we have $M$ is locally isometric to one of the following Lie groups with a left invariant metric: $SU(2)$, $SL(2,R)$, $E(2)$, $E(1,1)$.
%Moreover, we can find that the Ricci soliton is shrinking, steady, expanding for $SU(2)$($0<\mu<1$), $SL(2,R)$($\mu>1$), $E(2)$($\mu=1$), respectively. (A locally flat manifold is included as a special case in $E(2)$).
\vskip.15cm
Then, from now we find their associated potential vector fields $v$ explicitly.
For the non-Sasakian $(\alpha,\beta)$-manifold, we already know
$a=-\beta/2$ (constant) (Lemma 4.1 in \cite{BKP}) and the Ricci operator $S$ is given by
\begin{equation}\label{S}
S=-\beta I+\beta h+(2\alpha+\beta)\eta\otimes\xi,
\end{equation}
where $I$ denotes the identity transformation (Remark 3.2 in \cite{BKP}).
Thus, from (\ref{2.7}), (\ref{2.8}), (\ref{2.6}), (\ref{2}), (\ref{2.11}), and (\ref{S}) we have
\begin{equation}\label{2.12}
\left\{\begin{array}{lll}
e_1(f_2)+e_2(f_1)=0,\\
e_1(f_1)-\delta_1=0,\\
e_2(f_2)-\delta_2=0,\\
\xi(f_1)+\delta_3 f_2=0, \\
\xi(f_2)+\delta_4 f_1=0,
\end{array}\right.
\end{equation}
where $\delta_1=\beta-\beta\mu+2-2\mu^2$, $\delta_2=\beta+\beta\mu+2-2\mu^2$, $\delta_3=\beta/2+1+\mu$, and $\delta_4=-\beta/2+\mu-1$.

Moreover, we also observe that $e_1,e_2,\xi$ are geodesic vector fields.
 We adapt a normal coordinate system $(u^1,u^2,t)$ at $p\in M$, i.e.  $\exp_p \big(\sum_{i=1}^2  u^i(q)e_i+t(q)\xi\big)=q$.
From the last two equations of (\ref{2.12}), we establish 2nd order differential equations of constant coefficients:
\begin{equation}\label{2.20}
\frac{\partial^2 f_i}{\partial t^2}-\delta f_i=0
\end{equation}
at $p$, for $i=1,2$, where we have put $\delta=\delta_3\delta_4$.
Then we obtain their suitable solutions for $\delta>0,\ \delta<0,\ \delta=0$, respectively.
\vskip.2cm
\noindent
(I) For $\delta> 0$, we first get a general solution form of (\ref{2.20}):

\begin{equation*}\label{2.13}
\left\{\begin{array}{lll}
f_1(u^1,u^2,t)=A_1(u^1,u^2)\exp(\sqrt{\delta}t)+B_1(u^1,u^2)\exp(-\sqrt{\delta}t),\\
f_2(u^1,u^2,t)=A_2(u^1,u^2)\exp(\sqrt{\delta}t)+B_2(u^1,u^2)\exp(-\sqrt{\delta}t),
\end{array}\right.
\end{equation*}
where $A_i(u^1,u^2),B_i(u^1,u^2),$ $i=1,2,$ are smooth functions for $u^1,u^2$.
Further, from the first three equations of (\ref{2.12}) we may take
\begin{equation}\label{2.15}
\left\{\begin{array}{lll}
A_1(u^1,u^2)=&(\delta _1/2) u^1+(\delta _1/2) u^2+\tilde{C_1},\\ B_1(u^1,u^2)=&(\delta _1/2)u^1+(\delta _1/2) u^2+\tilde{D_1},\\
A_2(u^1,u^2)=&(-\delta _1/2) u^1+(\delta _2/2) u^2+\tilde{C_2},\\ B_2(u^1,u^2)=&(-\delta _1/2)u^1+(\delta _2/2) u^2+\tilde{D_2}.
\end{array}\right.
\end{equation}
Finally, reflecting the last two equations of (\ref{2.12}) again,
we obtain
\begin{equation}\label{2.17}
\left\{\begin{array}{lll}
 \tilde{C_1}=-\sqrt{\delta _3}C, \tilde{D_1}=\sqrt{\delta _3}D,\ \tilde{C_2}=\sqrt{\delta _4}C,\  \tilde{D_2}=\sqrt{\delta _4}D, \hbox{\ for } \delta _3>0, \delta _4>0,\\
 \tilde{C_1}=\sqrt{|\delta _3|}C, \tilde{D_1}=-\sqrt{|\delta _3|}D,\ \tilde{C_2}=\sqrt{|\delta _4|}C,\  \tilde{D_2}=\sqrt{|\delta _4|}D, \hbox{\ for } \delta _3<0, \delta _4<0.
\end{array}\right.
\end{equation}
where $C, D$ are arbitrary constants.
In similar ways, we have
\vskip.15cm
\noindent
(II) For $\delta<0$;
\begin{equation*}\label{2.14}
\left\{\begin{array}{lll}
f_1(u^1,u^2,t)=A_1(u^1,u^2)\cos(\sqrt{-\delta}t)+B_1(u^1,u^2)\sin(\sqrt{-\delta}t),\\
f_2(u^1,u^2,t)=A_2(u^1,u^2)\cos(\sqrt{-\delta}t)+B_2(u^1,u^2)\sin(\sqrt{-\delta}t).
\end{array}\right.
\end{equation*}
Here,
\begin{equation}\label{2.16}
\left\{\begin{array}{lll}
A_1(u^1,u^2)=&\delta _1 u^1+\delta _1 u^2+\tilde{C_1},\\
A_2(u^1,u^2)=&-\delta _1 u^1+\delta _2 u^2+\tilde{C_2},\\
B_1(u^1,u^2)=&\delta _1 u^1+\delta _1 u^2+\tilde{D_1}\\
B_2(u^1,u^2)=&-\delta _1 u^1+\delta _2 u^2 +\tilde{D_2},
\
\end{array}\right.
\end{equation}
and
\begin{equation}\label{2.18}
\left\{\begin{array}{lll}
 \tilde{C_1}=-\sqrt{\delta _3}C, \tilde{D_1}=\sqrt{\delta _3}D,\ \tilde{C_2}=-\sqrt{|\delta_4|}D,\ \tilde{D_2}=-\sqrt{|\delta _4|}C, \hbox{\ for } \delta _3>0, \delta _4<0,\\
 \tilde{C_1}=\sqrt{|\delta _3|}C, \tilde{D_1}=\sqrt{|\delta _3|}D,\ \tilde{C_2}=\sqrt{\delta_4}D,\ \tilde{D_2}=-\sqrt{\delta _4}C, \hbox{\ for } \delta _3<0, \delta _4>0,
\end{array}\right.
\end{equation}
where $C, D$ are arbitrary constants.
\vskip.15cm
\noindent
(III) For $\delta=0$;
\begin{equation*}\label{}
\left\{\begin{array}{lll}
f_1(u^1,u^2,t)=A_1(u^1,u^2)+B_1(u^1,u^2)t,\\
f_2(u^1,u^2,t)=A_2(u^1,u^2)+B_2(u^1,u^2)t,
\end{array}\right.
\end{equation*}
where
\begin{equation}\label{2.19}
\left\{\begin{array}{lll}
A_1(u^1,u^2)=&\delta _1 u^1+C,\\
A_2(u^1,u^2)=&\delta _2 u^2+D,\\
B_1(u^1,u^2)=&-\delta _3D, \\
B_2(u^1,u^2)=&-\delta _4C,
\end{array}\right.
\end{equation}
and $C,D$ are arbitrary constants. (For a locally flat manifold ($\mu=1$, $\beta=0$), we find $f_1=C-2D t$ and $f_2=D$.)
\vskip.3cm
Then, from the relations (\ref{2.15}) and (\ref{2.17}) for $\delta>0$ ((\ref{2.16}) and (\ref{2.18}) for $\delta<0$, (\ref{2.19}) for $\delta=0$, respectively) we can actually show that $bf_1+cf_2=0$ implies $b=c=0$ at $p$.
For example, for the case $\delta_3>0, \delta_4>0$ we compute $bf_1+cf_2=0$ at $p$:
$-b(p)\sqrt{\delta_3}(C-D)+c(p)\sqrt{\delta_4}(C+D)=0$ for arbitrary constants $C$ and $D$. Then, it implies at once $b=c=0$ since $\delta(=\delta_3\delta_4)\neq0.$
After all, $v=f_1e_1+f_2e_2$ are the desired potential vector fields.
This completes the proof of Theorem 2.
\vskip.3cm
We describe the Ricci solitons for the special cases $\beta=0$ by the table below.
\vskip.2cm
\noindent
\textrm{\red The special cases ($\beta=0$)}.
\vskip.2cm
\begin{tabular}{|r|l|l|l|}\hline
       & $0<\mu<1$ & $\mu>1$  & $\mu=1$  \\ \hline
 $\delta$ & negative & positive & zero \\ \hline
Potential vector field type & II & I & III \\ \hline
Ricci soliton & shirnking & expanding & steady \\ \hline
Lie group & SU(2) & SL(2,R) & E(2) \\ \hline
\end{tabular}

\section{\red The Sasakian case}
In the remainder of this paper, we treat a Ricci soliton for Sasakian 3-manifolds.
Geiges \cite{Gieges} proved the following result.
\vskip.15cm
\noindent
\begin{thm}
 A closed 3-manifold admits a Sasakian structure  if and only if it is diffeomorphic to a quotient of one of the spaces $S^3$, $\widetilde{SL(2,R)}$ or Nil by a discrete group of fixed point free isometries.
\end{thm}

 Let $M$ be a 3-dimensional unimodular Lie group with a
left-invariant Riemannian metric $g=\langle\cdot,\cdot\rangle$.
Then, according to a result of Perrone \cite{Perrone} $M$ admits its compatible left invariant Sasakian
structure if and only if there exists an orthonormal basis $\{e_1,
e_2, e_3\}$ of the Lie algebra $\mathfrak{m}$ such that

$$
[e_1,e_2]=2e_3,\ [e_2,e_3]=c_1e_1,\  [e_3, e_1]=c_1e_2.
$$
The Reeb vector filed $\xi$ is obtained
by left translation of $e_3$.
The contact distribution ${D}$ is
spanned by $e_1$ and $e_2$. It includes $SU(2)$, $SL(2,R)$, Nil(Heisenberg group), for $c_1$ positive, negative, zero, respectively.

By the Koszul formula, one can calculate the Levi-Civita
connection $\nabla$ in terms of the basis $\{e_1,e_2,e_3=\xi\}$ as
follows:

\begin{equation}
\begin{split}\label{4.0}
&\nabla_{e_1}e_1=0,\quad\nabla_{e_1}e_2=e_3,\quad
 \nabla_{e_1}e_3=-e_2,\\
&\nabla_{e_2}e_1=-e_3,\quad \nabla_{e_2}e_2=0,\quad \nabla_{e_2}e_3=e_1,\\
&\nabla_{e_3}e_1=(c_1-1)e_2,\quad
 \nabla_{e_3}e_2=-(c_1-1)e_1,\quad \nabla_{e_3}e_3=0.
\end{split}\end{equation}
From the above data, we are already aware that $e_1,e_2,e_3$ are geodesic vector fields. Moreover, we find that the basis  $\{e_1,e_2,e_3=\xi\}$ diagonalizes the Ricci operator, and we have
\begin{equation}\label{4.1}
\operatorname{Ric}(e_1,e_1)=\operatorname{Ric}(e_2,e_2)=2c_1-2,\ \
\operatorname{Ric}(\xi,\xi)=2.
\end{equation}
Namely, the Ricci operator $S$ is represented by $S=(2c_1-2)I+(4-2c_1)\eta\otimes\xi$, which says that $M$ is \textit{$\eta$-Einstein}.
We suppose that $M$ admits a  transversal Ricci soliton $(g,v,\lambda)$ with $v=f_1 e_1+f_2 e_2$, where $f_1$ and $f_2$ are smooth functions.
Then we have
\begin{equation}\label{Ricci111}
\frac{1}{2}\Big(g(\nabla_X v,Y)+g(\nabla_Y v,X)\Big) +  \operatorname{Ric}(X,Y) - \lambda g(X,Y) = 0.
\end{equation}
By similar computations in section 3, using (\ref{4.0}), (\ref{4.1}) and (\ref{Ricci111}), we finally obtain $\lambda=2$ (a shrinking Ricci soliton) and

\begin{equation}\label{4.2}
\left\{\begin{array}{lll}
e_1(f_2)+e_2(f_1)=0,\\
e_1(f_1)=e_2(f_2)=4-2c_1,\\
\xi(f_1)-(c_1-2) f_2=0, \\
\xi(f_2)+(c_1-2) f_1=0.
\end{array}\right.
\end{equation}
From the last two equations of (\ref{4.2}) we have
$$
\frac{\partial^2 f_i}{\partial t^2}+(c_1-2)^2 f_i=0
$$
for $i=1,2$.
Then, adapting a normal coordinate system $(u^1,u^2,t)$ associated with $\{e_1,
e_2, \xi\}$ and using similar arguments in section 3, we have a special solution for $c_1\neq 2$:
\begin{equation*}\label{4.3}
\left\{\begin{array}{lll}
f_1(u^1,u^2,t)=A_1(u^1,u^2)\cos(|c_1-2|t)+B_1(u^1,u^2)\sin(|c_1-2|t),\\
f_2(u^1,u^2,t)=A_2(u^1,u^2)\cos(|c_1-2|t)+B_2(u^1,u^2)\sin(|c_1-2|t).
\end{array}\right.
\end{equation*}
Here,
\begin{equation*}\label{4.5}
\left\{\begin{array}{lll}
A_1(u^1,u^2)=&(4-2c_1)(u^1+u^2)-D,\\
A_2(u^1,u^2)=&(4-2c_1)(-u^1+u^2)+C,\\
B_1(u^1,u^2)=&(4-2c_1)(u^1+u^2)+C,\\
B_2(u^1,u^2)=&(4-2c_1)(-u^1+u^2)+D,
\end{array}\right.
\end{equation*}
for $c_1>2$,
or
\begin{equation*}\label{4.4}
\left\{\begin{array}{lll}
A_1(u^1,u^2)=&(4-2c_1)(u^1+u^2)-D,\\
A_2(u^1,u^2)=&(4-2c_1)(-u^1+u^2)-C,\\
B_1(u^1,u^2)=&(4-2c_1)(u^1+u^2)+C,\\
B_2(u^1,u^2)=&(4-2c_1)(-u^1+u^2)-D,
\end{array}\right.
\end{equation*}
for $c_1<2$, where $C,D$ are arbitrary constants.
For the case $c_1=2$ ($M$ the unit sphere), we get a special $v=Ce_1+De_2$, where $C,D$ are arbitrary constants.
\vskip.2cm

For the Heisenberg group, we write down the Ricci soliton $(g,v,\lambda)$ in more detail.
\vskip.2cm
\noindent
\textrm{\red Example (Heisenberg group).}
For the Heisenberg group
\begin{equation*}
\mathbb{H}=\left\{
\begin{pmatrix}
1 & u^2 & t \\
0 & 1 & u^1 \\
0 & 0 & 1
\end{pmatrix} \vline u^1,u^2,t \in R\right\},
\end{equation*}
the contact form is $\eta=1/2(dt-u^2 du^1)$ and the Reeb vector field is $\xi=2\frac{\partial}{\partial t}$. The Riemannian metric $ds^2=\eta\otimes\eta+1/4((du^1)^2+(du^2)^2)$ is a left invariant metric and it is a
Sasakian metric associated with $\eta$ (cf. Example 4.5.1 in \cite{Blair}).
For the Heisenberg group $\mathbb{H}$ ($c_1=0$), we get a potential vector field $v=f_1e_1+f_2 e_2$, where
\begin{equation*}\label{4.5}
\left\{\begin{array}{lll}
f_1(u^1,u^2,t)=\big(4(u^1+u^2)-D\big)\cos(2t)+\big(4(u^1+u^2)+C\big)\sin(2t),\\
f_2(u^1,u^2,t)=\big(4(-u^1+u^2)-C\big)\cos(2t)+\big(4(-u^1+u^2)-D\big)\sin(2t),
\end{array}\right.
\end{equation*}
for arbitrary constants $C,D$.
\vskip.2cm

We close this paper by mentioning the Eight Geometries.
\vskip.2cm
\noindent
\textrm{\red Remark (Thurston geometry).} The eight model spaces of Thurston geometry
 are
©÷ 3-dimensional space forms $S^3, R^3, H^3$;
©÷ product spaces $S^2 \times R$ and $H^2 \times R$;
©÷ Nil and $\widetilde{SL(2,R)}$;
©÷ Sol.
(Class
$E(1.1)$ corresponds to Sol and $E(2)$ belongs to the flat geometry $R^3$.)
\vskip.5cm
\noindent
\textrm{\red Acknowledgement.}
This research was supported by Basic Science Research Program through the National Research
Foundation of Korea(NRF) funded by the Ministry of Education, Science and Technology
(2009-0071643).

\bigskip

\noindent Jong Taek Cho\\
Department of Mathematics,\\ Chonnam
National University,\\
Gwangju
500-757, Korea\\
E-mail address: {\tt  jtcho@chonnam.ac.kr}

\end{document}